\documentclass[psamsfonts]{tran-l}
\usepackage{amssymb,latexsym,graphicx}

\def\ds{\displaystyle}
\def\dist{\textrm{dist}}
\def\d{{\rm d}}
\def\e{{\rm e}}
\def\II{{\rm II}}

\def \R{\mathbb{R}}
\def \N{\mathbb{N}}
\def \Z{\mathbb{Z}}

\def \reff#1{(\ref{#1})}
\def \v#1{\vec{#1}}

\numberwithin{equation}{section}
\newtheorem{theorem}{Theorem}[section]
\newtheorem{lemma}[theorem]{Lemma}
\newtheorem{Remark}[theorem]{Remark}
\newenvironment{remark}{\begin{Remark}\rm}{\end{Remark}}
\newenvironment{varproof}
    {\rm \trivlist \item[\hskip \labelsep{\bf Proof}]}
    {\hspace*{\fill}$\Box$\endtrivlist}

\begin{document}

\title[Asymptotic zero distribution of multiple orthogonal polynomials]{Asymptotic
        zero distribution for a class of multiple orthogonal polynomials}
\thanks{This work was supported by INTAS project 03-51-6637, by FWO
        projects G.0455.04 and G.0184.02 and by OT/04/21 of K.U.Leuven}

\author{E.~Coussement, J.~Coussement, W.~Van Assche}
\address{Department of Mathematics\\ Katholieke Universiteit
        Leuven\\ Celestijnenlaan 200~B\\ 3001 Leuven, Belgium}
\thanks{The second author is a postdoctoral researcher at the K.U.Leuven (Belgium)}
\email{jonathan.coussement@wis.kuleuven.be}
\email{walter@wis.kuleuven.be}

\subjclass[2000]{Primary  33C45, 42C05; Secondary 15A18}

\keywords{Multiple orthogonal polynomials, Asymptotics}


\begin{abstract}
We establish the asymptotic zero distribution for polynomials
generated by a four-term recurrence relation with varying
recurrence coefficients having a particular limiting behavior. The
proof is based on ratio asymptotics for these polynomials. We can
apply this result to three examples of multiple orthogonal
polynomials, in particular Jacobi-Pi\~neiro, Laguerre~I and the
example associated with Macdonald functions. We also discuss an
application to Toeplitz matrices.
\end{abstract}

\maketitle

\section{Introduction}
\label{intro}

Let $\mu$ be a positive measure on the real line for which the
support is not finite and all the moments exist. The corresponding
monic orthogonal polynomial $P_n$ of degree $n$ is then defined by
\begin{equation}
\int x^m P_n(x)\, \d \mu(x)=0,\qquad k=0,\ldots,n-1,
\end{equation}
with $P_0\equiv 1$ and $P_{-1}\equiv 0$.  A well-known fact is
that such polynomials satisfy a three-term recurrence relation of
the form
\begin{equation}
\label{recursie OP} %
zP_n(z)=P_{n+1}(z)+b_n P_n(z)+a_n^2 P_{n-1}(z),\qquad a_n>0,\
b_n\in\R,
\end{equation}
with initial conditions $P_0\equiv 1$ and $P_{-1}\equiv 0$.

An object of frequent study is the {\em asymptotic zero
distribution} of the zeros for a sequence of orthogonal
polynomials. The zeros of the polynomials $P_n$, generated by
\reff{recursie OP}, are real and simple \cite{Chihara}. With each
polynomial $P_n$ we can associate the normalized zero counting
measure
\begin{equation}
\label{nu} %
\nu (P_n):=\frac{1}{n}\sum_{P_n(x)=0} \delta_{x},
\end{equation}
where $\delta_x$ is the Dirac point mass at $x$.  If $\lim_{n\to
\infty}\nu (P_n)= \nu$, by which we mean that
\[\lim_{n\to \infty}\int f\, \d\nu (P_n)=\int f\, \d\nu\]
for every bounded and continuous function $f$ on $\mathbb{R}$
({\em weak-$\star$ convergence}), then we call the probability
measure $\nu$ the asymptotic zero distribution of the sequence
$\{P_n\}_{n\ge 0}$. One of the famous results in this context is
the following.
\begin{theorem}[see, e.g., \cite{Nevai,Walter}]
\label{non varying OP} %
Suppose that the recurrence coefficients $a_n>0$ and $b_n\in\R$
have the limits $a>0$ and $b\in\R$, respectively. The polynomials
$P_n$, generated by \reff{recursie OP}, then have the asymptotic
zero distribution $\omega_{[\gamma,\delta]}$ with density
\begin{equation}
\label{arcsine} %
\frac{\d \omega_{[\gamma,\delta]}}{\d x}(x)= \left\{ %
\begin{array}{ll}
{\ds \frac{1}{\pi \sqrt{(\delta-x)(x-\gamma)}}}, & \qquad x \in
[\gamma,\delta],\\[2ex]
0, & \qquad \mbox{elsewhere},
\end{array} \right.
\end{equation}
where $\gamma=b-2a$ and $\delta=b+2a$.
\end{theorem}
\begin{remark}
The measure $\omega_{[\gamma,\delta]}$ is known as the arcsine
measure on $[\gamma,\delta]$. It also minimizes the logarithmic
energy of the interval $[\gamma,\delta]$ \cite{Saff}.
\end{remark}
Recently,  the result in Theorem~~\ref{non varying OP} was
extended to the case of varying recurrence coefficients. Here the
notation $\lim_{n/N \to t} Y_{n,N}= Y$ denotes the property that
in the doubly indexed sequence $Y_{n,N}$ we have $\lim_{j \to
\infty} Y_{n_j,N_j} = Y$ whenever $n_j$ and $N_j$ are two
sequences of natural numbers such that $N_j \to \infty$ and
$n_j/N_j \to t$ as $j \to \infty$.
\begin{theorem}[Kuijlaars, Van Assche \cite{Kuijlaars1}]
\label{varying OP} %
Let for each $N\in \N$, two sequences
$\{a_{n,N}\}_{n=1}^\infty$, $a_{n,N}>0$, and
$\{b_{n,N}\}_{n=0}^\infty$, $b_{n,N}\in \R$, of recurrence
coefficients be given. Furthermore, suppose there exist two
continuous functions $a:(0,+\infty)\to [0,+\infty)$,
$b:(0,+\infty)\to \R$, such that
\begin{equation}
\lim_{n/N \to t} a_{n,N}=a(t),\qquad \lim_{n/N \to t}
b_{n,N}=b(t),\qquad t>0,
\end{equation}
and define $\gamma(t):=b(t)-2a(t)$, $\delta:=b(t)+2a(t)$, $t>0$.
For the (orthogonal) polynomials generated by the recurrence
\begin{equation}
zP_{n,N}(z)=P_{n+1,N}(z)+b_{n,N} P_{n,N}(z)+a_{n,N}^2
P_{n-1,N}(z),
\end{equation}
with initial conditions $P_{0,N}\equiv 1$ and $P_{-1,N}\equiv 0$,
we then have
\begin{equation}
\lim_{n/N \to t} \nu(P_{n,N})=\frac{1}{t}\int_0^t
\omega_{[\gamma(s),\delta(s)]}\,\d s,\qquad t>0.
\end{equation}
Here $\omega_{[\gamma,\delta]}$ is defined by \reff{arcsine} if
$\gamma<\delta$ and by $\delta_\gamma$ if $\gamma=\delta$.
\end{theorem}
\begin{remark}
More recently, Theorem~\ref{varying OP} was generalized to
measurable functions $a$ and $b$ \cite{Kuijlaars3}.
\end{remark}


In this paper we present a (conditional) theorem giving the
asymptotic zero distribution for polynomials satisfying a
four-term recurrence relation of the form
\begin{equation}
\label{recursie intro} %
z P_{n,N}(z) = P_{n+1,N}(z) + b_{n,N}
P_{n,N}(z) + c_{n,N} P_{n-1,N}(z) + d_{n,N} P_{n-2,N}(z),
\end{equation}
where the varying recurrence coefficients have some particular
limiting behavior.  So, in a sense it extends Theorem~\ref{varying
OP}. Such a four-term recurrence relation appears in the theory of
multiple orthogonal polynomials of Type~$\II$. These are a
generalization of orthogonal polynomials which arises naturally in
Hermite-Pad\'e approximation of a system of (Markov) functions
\cite{Bruin1,Bruin2,Mahler}. In particular, they satisfy
orthogonality conditions with respect to several positive measures
\cite{Aptekarev2,Nikishin,Els1}. Some of their applications are
situated in diophantine number theory, rational approximation,
spectral and scattering problems for higher-order difference
equations and some associated dynamical systems, see, e.g.,
\cite{Aptekarev3,bultheel,Ismail,Walter3}. Recently they also
appeared in random matrix theory for matrix ensembles with
external source \cite{bleher3,bleher1,bleher2} and Wishard
ensembles \cite{bleher4}. The particular limiting behavior which
we are considering appears in the examples Jacobi-Pi\~neiro,
Laguerre~I \cite{Els1} and the example associated with Macdonald
functions \cite{Walter2}.

In Subsection~\ref{main theorem} we state our main theorem. Next,
in Subsection~\ref{MOP} and Subsection~\ref{TM} we apply this
result to the examples of multiple orthogonal polynomials
mentioned above and some particular kind of Toeplitz matrices.  In
Section~\ref{RA} we discuss a theorem on ratio asymptotics for
monic polynomials satisfying the recurrence \reff{recursie intro}.
This will help us to prove our main theorem in Section~\ref{proof
main theorem}.

\section{Statement of results}

\subsection{Main theorem}
\label{main theorem}

We will study doubly indexed sequences of polynomials
$\{P_{n,N}\}$, generated by a four-term recurrence of the form
\begin{equation}
\label{four term recursie} %
z P_{n,N}(z) = P_{n+1,N}(z) + b_{n,N} P_{n,N}(z) + c_{n,N}
P_{n-1,N}(z) + d_{n,N} P_{n-2,N}(z),
\end{equation}
with the initial conditions $P_{0,N}\equiv 1$, $P_{-1,N}\equiv 0$
and $P_{-2,N}\equiv 0$ and real recurrence coefficients. In
particular, our main theorem gives an explicit expression for the
asymptotic zero distributions
    \[\lim_{n/N\to t}\nu (P_{n,N}),\qquad t>0,\]
with some conditions on the zeros of the $P_ {n,N}$ and some
particular limiting behavior for the recurrence coefficients. As
mentioned in the introduction the limit is taken over any sequence
$\{\nu (P_{n_j,N_j})\}_{j\ge 1}$ for which $n_j\to \infty$,
$N_j\to \infty$ and $n_j/N_j\to x$ as $j\to \infty$. We will use
this notation throughout the rest of this paper.
\begin{theorem}
\label{as zero distr} %
Let for each $N \in \mathbb{N}$ three sequences
$\{b_{n,N}\}_{n=0}^{\infty}$, $\{c_{n,N}\}_{n=1}^{\infty}$ and
$\{d_{n,N}\}_{n=2}^{\infty}$ of real recurrence coefficients be
given and assume that there exists a continuous function
$\alpha:[0,+\infty)\to [0,+\infty)$ such that, for $t>0$,
\begin{equation}
\label{limbcd met t} %
\lim_{n/N \to t} b_{n,N} = 3\beta(t), %
\quad \lim_{n/N\to t} c_{n,N} = 3\beta(t)^2, %
\quad \lim_{n/N \to t} d_{n,N} = \beta(t)^3,
\end{equation}
with $\beta(t)=\frac{4\alpha(t)}{27}$.
Let $P_{n,N}$ be the monic polynomials generated by the recurrence
\reff{four term recursie} and suppose these polynomials $P_{n,N}$
have real simple zeros $x_{1}^{n,N}< \ldots < x_{n}^{n,N}$
satisfying the interlacing property $ x_{j}^{n+1,N} < x_{j}^{n,N}
< x_{j+1}^{n+1,N}$, for all $n,N \in \mathbb{N}$, $j=1, \ldots,
n$. Then
\begin{equation}
\label{limitnu} %
\lim_{n/N\to t}\nu (P_{n,N})=\frac{1}{t}\int_0^t
\upsilon_{[0,\alpha(s)]}\,\d s,\qquad t>0,
\end{equation}
where $\upsilon_{[0,\alpha]}$ is defined by $\delta_0$ if
$\alpha=0$ and $\frac{\d \upsilon_{[0,\alpha]}}{\d
x}(x)=\frac{1}{\alpha}\frac{\d \upsilon_{[0,1]}}{\d
x}(\frac{x}{\alpha})$, with
\begin{equation}
\label{upsilon} %
\frac{\d \upsilon_{[0,1]}}{\d x}(x)= \left\{ %
\begin{array}{ll}
{\ds \frac{\sqrt{3}}{4 \pi}\,\frac{(1+\sqrt{1-x})^{1/3}+
(1-\sqrt{1-x})^{1/3}}{x^{2/3} \sqrt{1-x}}}, & \qquad x \in
(0,1),\\[2ex]
0, & \qquad \mbox{elsewhere},
\end{array} \right.
\end{equation}
if $\alpha>0$.
\end{theorem}
\begin{remark}
\label{nu t} %
Denote by $\nu_t$ the right hand side of \reff{limitnu}. This
measure acts on an arbitrary Borel set $E$ like
\[\nu_t(E)=\frac{1}{t} \int_0^t\upsilon_{[0,\alpha(s)]}(E)\, \d s,
\qquad E\in \mathcal{B}(\mathbb{R}).\]
Now suppose that for each $x\ge 0$ the set
$\{s\ge 0 \,|\, \ x\le \alpha(s)\}$
is an interval, which we denote by $[t_-(x),t_+(x)]$. The density
of the measure $\nu_t$ is then
\begin{equation}
\frac{\d \nu_t}{\d
x}(x)=\frac{1}{t}\int_{\min(t,t_-(x))}^{\min(t,t_+(x))} \frac{\d
\upsilon_{[0,\alpha(s)]}}{\d x}(x) \,\d s.
\end{equation}
This will be the case in each of the examples we present in this
paper.
\end{remark}

\begin{remark}
Comparing Theorem~\ref{as zero distr} with Theorem~\ref{varying
OP} we see that the measure $\upsilon_{[0,1]}$ plays the role of
the arcsine measure in the case of orthogonal polynomials
(satisfying a three-term recurrence relation).  The density
\reff{upsilon} again has the behavior $c_1(1-x)^{-1/2}$ as $x
\uparrow 1$, but has a different behavior $c_2x^{-2/3}$ as
$x\downarrow 0$.
\end{remark}

\begin{remark}
The measure $\upsilon_{[0,1]}$ coincides (after a cubic
transformation) with the asymptotic zero distribution of Faber
polynomials associated with the 3-cusped hypocycloid
\cite{Kuijlaars2}.
\end{remark}

\begin{figure}[t]
\begin{center}
\includegraphics[scale=0.4]{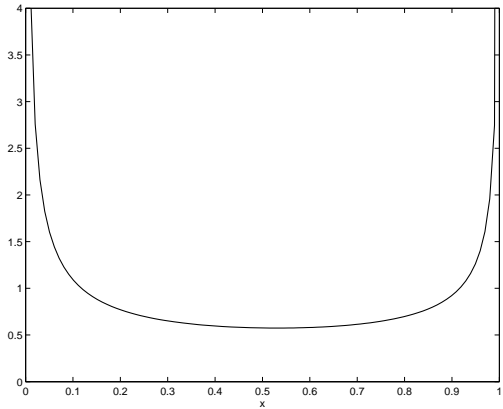}
\end{center}
\caption{\label{maat v}{\em The density of the measure
$\upsilon_{[0,1]}$.}}
\end{figure}

\subsection{Application to multiple orthogonal polynomials}
\label{MOP}

There are two types of multiple orthogonal polynomials but we will
only consider type~$\II$. Let $\mu_1,\ldots,\mu_r$, $r\in \N$, be
a set of positive measures on the real line for which the support
is not finite and all the moments exist. Furthermore, let $\v n =
(n_1,n_2,\ldots,n_r)$ be a vector of $r$ nonnegative integers,
which is a {\em multi-index} with length $|\v n| :=
n_1+n_2+\cdots+n_r$. A multiple orthogonal polynomial $P_{\v n}$
of type II with respect to the multi-index $\v n$, is a
(nontrivial) polynomial of degree $\leq |\v n|$ which satisfies
the orthogonality conditions
\begin{equation}
\label{stelseltypeII} %
\int x^m P_{\v n}(x) \,\d\mu_j (x) =0, \qquad 0\le m\le
n_j-1,\qquad j=1,\ldots ,r.
\end{equation}
A basic requirement in the study of multiple orthogonal
polynomials is that the system (\ref{stelseltypeII}) has a unique
solution (up to a scalar multiplicative constant) of degree $|\v
n|$. We call $\v n$ a {\em normal index} for $\mu_1,\ldots,\mu_r$
if any solution of (\ref{stelseltypeII}) has exactly degree $|\v
n|$ (which implies uniqueness). If all the multi-indices are
normal then the system of measures is called {\em perfect}. Some
famous classes of perfect systems are the Angelesco systems,
Nikishin systems ( for $r=2$) and AT~systems, see, e.g.,
\cite{Nikishin,Els1}.

Multiple orthogonal polynomials of type~$\II$ satisfy a recurrence
relation of order $r+1$. In particular, if we set $r=2$ and
consider {\em proper multi-indices} $\v \nu_n=(m+s,m)$, $n\in\N
\cup \{0\}$, where $n=2m+s$, $s\in \{0,1\}$, then the polynomials
$P_n:=P_{\v \nu_n}$ satisfy a four-term recurrence relation of the
form
\begin{equation}
\label{four term recursie MOP} %
z P_{n}(z) = P_{n+1}(z) + b_{n} P_{n}(z) + c_{n} P_{n-1}(z) +
d_{n} P_{n-2}(z),
\end{equation}
with the initial conditions $P_{0}\equiv 1$, $P_{-1}\equiv 0$ and
$P_{-2}\equiv 0$. For three examples known in the literature the
recurrence coefficients in \reff{four term recursie MOP} have the
particular limiting behavior \reff{limbcd met t}, possibly after
some re-scaling.  In each of these examples the measures form an
AT~system on an interval $\Delta\subseteq \R$. It is then known
that the zeros of the polynomials $P_n$ are simple, lie in
$\Delta$ \cite{Nikishin,Els1} and satisfy the interlacing property
\cite{Aptekarev1}.  So, it is possible to apply Theorem~\ref{as
zero distr}.

\subsubsection{Jacobi-Pi\~neiro}
The Jacobi-Pi\~neiro polynomials are the multiple orthogonal
polynomials for the system of orthogonality measures \[\d
\mu_j(x)=x^{\alpha_j}(1-x)^\beta\,\d x, \quad j=1,2,\] on the
interval $[0,1]$ with $\alpha_1,\alpha_2,\beta>-1$ and
$\alpha_2-\alpha_1\notin \Z$. In \cite{Els1} it was shown that the
monic Jacobi-Pi\~neiro polynomials with respect to proper
multi-indices, which we denote by $P_n^{\alpha_1,\alpha_2;\beta}$,
satisfy a recurrence relation of the form \reff{four term recursie
MOP} for which
\begin{equation}
\lim_{n\to \infty}b_n=3\left(\frac{4}{27}\right),\qquad \lim_{n\to
\infty}c_n=3\left(\frac{4}{27}\right)^2,\qquad \lim_{n\to
\infty}d_n=\left(\frac{4}{27}\right)^3.
\end{equation}
By Theorem~\ref{as zero distr} with $\alpha(t)=1$, $t>0$, we then
easily obtain the following result.
\begin{theorem}
The Jacobi-Pi\~neiro polynomials $P_n^{\alpha_1,\alpha_2;\beta}$
have the asymptotic zero distribution $\upsilon_{[0,1]}$, defined
as in \reff{upsilon}.
\end{theorem}

\subsubsection{Multiple Laguerre I}

The multiple Laguerre polynomials of the first kind are orthogonal
with respect to the system of measures
\[\d \mu_j(x)=x^{\alpha_j}\e^{-x}\,\d x, \quad j=1,2,\] on
$[0,+\infty)$ with $\alpha_1,\alpha_2>-1$ and
$\alpha_2-\alpha_1\notin \Z$. Denote the monic multiple Laguerre~I
polynomials with respect to proper multi-indices by
$L_n^{\alpha_1,\alpha_2}$. These satisfy a four-term recurrence
relation of the form \reff{four term recursie MOP} where, for
$t>0$,
\begin{equation}
\lim_{n/N\to t}\frac{b_n}{N}=3\left(\frac{t}{2}\right),\qquad
\lim_{n/N\to t}\frac{c_n}{N^2}=3\left(\frac{t}{2}\right)^2,\qquad
\lim_{n/N\to t}\frac{d_n}{N^3}=\left(\frac{t}{2}\right)^3,
\end{equation}
see \cite{Els1}. The following theorem is then a corollary of
Theorem~\ref{as zero distr}.
\begin{theorem}
\label{th laguerre} %
For the multiple Laguerre polynomials of the
first kind the limit
\begin{equation}
\nu_t^L:=\lim_{n/N\to
t}\sum_{L_n^{\alpha_1,\alpha_2}(x)=0}\delta_{x/N},\qquad t>0,
\end{equation}
exists and has the density
\begin{equation}
\label{density L}
\frac{\d \nu_t^L}{\d x}(x)= \left\{ %
\begin{array}{ll}
\frac{8}{27 t} \, g\left(\frac{8x}{27t}\right), & \qquad x \in
(0,\frac{27 t}{8}),\\[2ex]
0, & \qquad \mbox{elsewhere},
\end{array} \right.
\end{equation}
where, for $y\in (0,1)$,
\[g(y)=\frac{3\sqrt{3}}{16\pi}
\frac{(1+3\sqrt{1-y})(1-\sqrt{1-y})^{1/3}-(1-3\sqrt{1-y})(1+\sqrt{1-y})^{1/3}}{y^{2/3}}.\]
\end{theorem}

\begin{proof}
If we define $\tilde L_n^{\alpha_1,\alpha_2}(z):=
L_n^{\alpha_1,\alpha_2}(Nz)/N^n$, then the polynomials satisfy a
recurrence relation of the form \reff{four term recursie} and the
asymptotic property \reff{limbcd met t} with
$\alpha(t)=\frac{27t}{8}$, $t>0$. So, applying Theorem~\ref{as
zero distr} we get
\[
\frac{\d \nu_t^L}{\d x}(x)  =  %
\frac{1}{t} \int\limits_{\frac{8x}{27}}^t \frac{8}{27s} \frac{\d
\upsilon_{[0,1]}}{\d x}\left(\frac{8x}{27s}\right)\d s
 =  %
\frac{1}{t}\frac{8}{27} \int\limits_{\frac{8x}{27t}}^1
\frac{1}{u}\frac{\d \upsilon_{[0,1]}}{\d x}\left(u\right)\d u.
\]
Set $x \in (0,\frac{27 t}{8})$ and $z=\sqrt{1-\frac{8x}{27t}}$.
Applying the substitution $y\leftrightarrow \sqrt{1-u}$ we then
obtain
\begin{eqnarray*}
\frac{\d \nu_t^L}{\d x}(x) & = & %
\frac{1}{t}\frac{4\sqrt{3}}{27\pi} \int_{0}^z \left((1+y)^{-4/3}(1-y)^{-5/3}+(1+y)^{-5/3}(1-y)^{-4/3}\right)\d y\\
& = & %
\frac{1}{t}\frac{4\sqrt{3}}{27\pi} \int_{-z}^z
(1+y)^{-5/3}(1-y)^{-4/3}\, \d y \\
& = & %
\left. \frac{1}{t}\frac{\sqrt{3}}{18\pi}\,
\frac{1+3y}{(1+y)^{2/3}(1-y)^{1/3}} \, \right|_{-z}^z.
\end{eqnarray*}
This completes the proof.
\end{proof}

\begin{figure}[t]
\begin{center}
\includegraphics[scale=0.4]{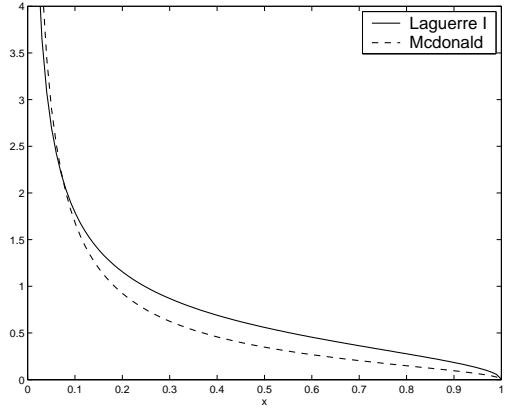}
\end{center}
\caption{\label{maat L en M}{\em The densities of the measures
$\nu_t^L$ and $\nu_t^M$, see \reff{density L} and \reff{density
M}, with $t=\frac{8}{27}$ and $t=\frac{2}{3\sqrt{3}}$
respectively.}}
\end{figure}

\subsubsection{Multiple orthogonal polynomials associated with Macdonald
functions}

In \cite{Walter2} one considered multiple orthogonal polynomials
with respect to the orthogonality measures
\[\d \mu_1(x)=x^\kappa\rho_\gamma(x)\,\d x, \quad \d \mu_2(x)=x^\kappa\rho_{\gamma+1}(x)\,\d x\]
on $(0,+\infty)$ with $\kappa>-1$, $\gamma\ge 0$ and
\[\rho_\gamma(x)=2x^{\gamma/2}K_\gamma(2\sqrt{x}),\qquad x>0,\]
where $K_\gamma$ is the modified Bessel of the second kind, also
known as the Macdonald function \cite[p.~374]{Abra}. For the type
II multiple polynomials with respect to proper multi-indices,
$P_n^{\gamma;\kappa}$, the recurrence coefficients in \reff{four
term recursie MOP} are known \cite[Theorem~4]{Walter2}. In
particular, for $t>0$,
\begin{equation}
\label{lim macdonald}
\lim_{n/N\to t}\frac{b_n}{N^2}=3t^2, %
\qquad \lim_{n/N\to t}\frac{c_n}{N^4}=3t^4, %
\qquad \lim_{n/N\to t}\frac{d_n}{N^6}=t^6.
\end{equation}
Theorem~\ref{as zero distr} then implies the following asymptotic
result for the zeros of these polynomials.
\begin{theorem}
For the multiple orthogonal polynomials associated with Macdonald
functions the limit
\begin{equation}
\nu_t^M:=\lim_{n/N\to
t}\sum_{P_n^{\gamma;\kappa}(x)=0}\delta_{x/N^2}
\end{equation}
exists and has the density
\begin{equation}
\label{density M}
\frac{\d \nu_t^M}{\d x}(x)= \left\{ %
\begin{array}{ll}
\frac{4}{27 t^2} \, h\left(\frac{4x}{27t^2}\right), & \qquad x
\in (0,\frac{27 t^2}{4}),\\[2ex]
0, & \qquad \mbox{elsewhere},
\end{array} \right.
\end{equation}
where
\[h(y)=\frac{3\sqrt{3}}{4\pi}\,
\frac{(1+\sqrt{1-y})^{1/3}-(1-\sqrt{1-y})^{1/3}}{y^{2/3}},
\qquad y\in (0,1).\]
\end{theorem}

\begin{proof}
Define the polynomials $\tilde P_n^{\gamma;\kappa}(z):=
P_n^{\gamma;\kappa}(N^2z)/N^{2n}$. By \reff{lim macdonald} these
satisfy a recurrence relation of the form \reff{four term
recursie} having the asymptotic property \reff{limbcd met t} with
$\alpha(t)=\frac{27t^2}{4}$, $t>0$. Set $x \in (0,\frac{27
t^2}{4})$ and $z=\frac{4x}{27 t^2}$. Applying Theorem~\ref{as zero
distr} then gives
\[
\frac{\d \nu_t^M}{\d x}(x)  =  %
\frac{1}{t} \int\limits_{\frac{2}{3}\sqrt{\frac{x}{3}}}^t
\frac{4}{27s^2} \frac{\d \upsilon_{[0,1]}}{\d
x}\left(\frac{4x}{27s^2}\right)\d s
 =  %
\frac{1}{t^2}\frac{2}{27\sqrt{z}} \int\limits_{z}^1
\frac{1}{\sqrt{u}}\frac{\d \upsilon_{[0,1]}}{\d x}(u)\, \d u.
\]
Similarly as in Theorem~\ref{th laguerre} we apply the
substitution $y\leftrightarrow \sqrt{1-u}$ and get
\begin{eqnarray*}
\frac{\d \nu_t^M}{\d x}(x) & = & %
\frac{1}{t^2}\frac{\sqrt{3}}{27\pi}\frac{1}{\sqrt{z}}
\int_{-\sqrt{1-z}}^{\sqrt{1-z}}(1+y)^{-5/6}(1-y)^{-7/6}\d y\\
& = & %
\left. \frac{1}{t^2}\frac{1}{3\sqrt{3}\pi}\frac{1}{\sqrt{z}}\,
\left(\frac{1+y}{1-y}\right)^{1/6} \,
\right|_{-\sqrt{1-z}}^{\sqrt{1-z}}.
\end{eqnarray*}
From this we easily obtain \reff{density M}.
\end{proof}

\subsection{Application to Toeplitz matrices}
\label{TM}

Set $\alpha>0$, $\beta=\frac{4\alpha}{27}$, and define the
Toeplitz matrices
\begin{eqnarray*}
T_n^\alpha&:=& \left(
\begin{array}{cccccc}
3\beta&1 & 0 & \ldots & \ldots &  0\\
3\beta^2&3\beta &1 &\ddots &  &\vdots \\
\beta^3&3\beta^2 &3\beta &1 & \ddots  & \vdots\\
0&\beta^3&3\beta^2 &3\beta  & \ddots & 0\\
\vdots&\ddots &\ddots &\ddots  &\ddots &1 \\
0&\ldots &0 &\beta^3 &3\beta^2 & 3\beta
\end{array}
\right)\in \mathbb{R}^{n\times n},\qquad n\in \N.
\end{eqnarray*}
Note that the eigenvalues of $T_n^\alpha$ coincide with the zeros
of the monic polynomial $Q_n^\alpha(z)=\det(zI_n-T_n^\alpha)$,
$n\in\mathbb{N}$. These polynomials satisfy the recurrence
relation
\begin{equation}
\label{recursie Q} %
z\, Q_n^\alpha(z) = Q_{n+1}^\alpha(z) + 3\beta\, Q_{n}^\alpha(z) +
3\beta^2\, Q_{n-1}^\alpha(z) + \beta^3\, Q_{n-2}^\alpha(z),
\end{equation}
with $Q_{0}^\alpha\equiv 1$ and $Q_{-1}^\alpha\equiv
Q_{-2}^\alpha\equiv 0$. The following asymptotic result for the
eigenvalues of the matrices $T_n^\alpha$ then follows from
Theorem~\ref{as zero distr}.
\begin{theorem}
\label{limit Toeplitz} %
The limiting eigenvalue distribution of the matrices $T_n^\alpha$,
with $\alpha>0$, is given by the measure $\upsilon_{[0,\alpha]}$,
defined as in Theorem~\ref{as zero distr}.
\end{theorem}

\begin{proof}
The homogeneous recurrence relation
\[0 = Q_{n+1}^\alpha(0) + 3\beta\,
Q_{n}^\alpha(0) + 3\beta^2\, Q_{n-1}^\alpha(0) + \beta^3\,
Q_{n-2}^\alpha(0)
\]
with $Q_{0}^\alpha(0)= 1$ and $Q_{-1}^\alpha(0)= Q_{-2}^\alpha(0)=
0$, has the solution
\[Q_n^\alpha(0)=(-\beta)^n\left(1+\frac{3n}{2}+\frac{n^2}{2}\right), \qquad n\in\mathbb{N}.\]
So, since $\beta=\frac{4\alpha}{27}>0$ all the $T_n^\alpha$ are
non-singular. Next, define
\begin{eqnarray*}
\tilde{T}_n^\alpha&=& \left(
\begin{array}{ccccccc}
1&0 & \ldots &\ldots & \ldots & \ldots &0 \\
3\beta&1 & \ddots & & & & \vdots\\
3\beta^2&3\beta &1 &\ddots & & &\vdots \\
\beta^3&3\beta^2 &3\beta &1 & \ddots & & \vdots\\
0&\beta^3&3\beta^2 &3\beta &1 & \ddots & \vdots\\
\vdots&\ddots &\ddots &\ddots &\ddots &\ddots &0 \\
0&\ldots &0 &\beta^3 &3\beta^2 & 3\beta&1
\end{array}
\right)\in \mathbb{R}^{(n+1)\times(n+1)},
\end{eqnarray*}
and notice that $\tilde{T}_n^\alpha=(\tilde{A}_n^\alpha)^3$ with
\begin{eqnarray*}
\tilde{A}_n^\alpha &=& \left(
\begin{array}{ccccc}
1&0& \ldots& \ldots &0\\ \beta&1&\ddots& &\vdots\\ 0&\beta&1& \ddots & \vdots\\
\vdots& \ddots& \ddots& \ddots& 0\\ 0 & \ldots&0 &\beta&1
\end{array}
\right)\in \mathbb{R}^{(n+1)\times(n+1)}.
\end{eqnarray*}
Since $\tilde{A}_n^\alpha$ is totally non-negative,
$\tilde{T}_n^\alpha$ is totally non-negative, see, e.g.,
\cite[p.~74, $1^\circ$]{Gant}, and so also $T_n^\alpha$. By
\cite[p.~100, Theorem~10]{Gant} we then get that the $T_n^\alpha$
are oscillation matrices.
Consequently, the zeros of the polynomials $Q_n$ are simple and
positive \cite[p.~87, Theorem~6]{Gant} and satisfy the interlacing
property \cite[p.~107, Theorem~14]{Gant}. The theorem then easily
follows from Theorem~\ref{as zero distr} and the recurrence
\reff{recursie Q}.
\end{proof}

\begin{remark}
The polynomials satisfying the recurrence relation with constant
coefficients \reff{recursie Q} are the multiple Chebyshev
polynomials of the second kind after a cubic transformation
\cite{Cheikh,Douak}. These are an example of multiple orthogonal
polynomials of type~$\II$ extending the well-known Chebyshev
polynomials of the second kind \cite{Chihara}. The corresponding
orthogonality measures can be found in
\cite[Corollary~4.2]{Cheikh}, \cite[Theorem~4.1~]{Douak}.
\end{remark}

\section{Ratio asymptotics}
\label{RA}

In \cite{Kuijlaars1} Kuijlaars and Van Assche have proven a
theorem that gives explicit ratio asymptotics for orthogonal
polynomials with converging varying recurrence coefficients. In
this section we give an extension of this result to polynomials
satisfying a (specific) four-term recurrence relation instead of a
three-term recurrence relation.

\begin{theorem}
\label{ratioas} %
Suppose we have for each $N \in \mathbb{N}$ sequences
$\{b_{n,N}\}_{n=0}^{\infty}$, $\{c_{n,N}\}_{n=1}^{\infty}$ and
$\{d_{n,N}\}_{n=2}^{\infty}$ of real recurrence coefficients and
let $P_{n,N}$ be the monic polynomials generated by the recurrence
\begin{equation}
\label{recursie ratio as} %
z P_{n,N}(z) = P_{n+1,N}(z) + b_{n,N}
P_{n,N}(z) + c_{n,N} P_{n-1,N}(z) + d_{n,N} P_{n-2,N}(z),
\end{equation}
with $P_{0,N}\equiv 1$, $P_{-1,N}\equiv 0$ and $P_{-2,N}\equiv 0$.
Assume that for some fixed $t>0$ the recurrence coefficients have
the limits
\begin{equation}
\label{limbcd} %
\lim_{n/N \to t} b_{n,N} = 3\left(\frac{4\alpha}{27}\right), %
\quad \lim_{n/N\to t} c_{n,N} = 3\left(\frac{4\alpha}{27}\right)^2, %
\quad \lim_{n/N \to t} d_{n,N} =
\left(\frac{4\alpha}{27}\right)^3,
\end{equation}
with $\alpha \ge 0$.
Furthermore, assume that the polynomials $P_{n,N}$ have real
simple zeros $x_{1}^{n,N}< \ldots < x_{n}^{n,N}$ satisfying the
interlacing property $ x_{j}^{n+1,N} < x_{j}^{n,N} <
x_{j+1}^{n+1,N}$, for all $n,N \in \mathbb{N}$, $j=1, \ldots, n$.
Moreover, suppose that for some $t^{*}>t$, there exist $m \le 0$,
$M \ge \alpha$ such that all zeros of $P_{n,N}$ belong to $[m,M]$
whenever $n \leq t^{*} N$.
Then
\begin{equation}
\label{ratiolimit} %
\lim_{n/N \to t} \frac{P_{n,N}(z)}{P_{n+1,N}(z)} =
\left\{
\begin{array}{ll}
{\ds \frac{1}{\alpha}\,\phi\left(\frac{z}{\alpha}\right)}, &
\qquad \alpha>0,\\[3ex]
{\ds \frac{1}{z}}, & \qquad \alpha=0,
\end{array}
\right.
\end{equation}
uniformly on compact subsets of $\mathbb{C}\setminus [m,M]$, where
$\phi$ is  defined by
\begin{equation}
\label{phidefinitie} %
\phi(z) := \frac{27}{4}\left(\frac{3\, \omega_3}{2}\, z^{1/3}
\left\{ \omega_3 (-1+ \sqrt{1-z})^{1/3} + (-1- \sqrt{1-z})^{1/3}
\right\}-1\right)
\end{equation}
with $\omega_3= e^{\frac{2 \pi i}{3}}$ and
\begin{eqnarray}
\nonumber
\sqrt{\rho e^{i \theta}} &=& \rho^{1/2}\, e^{\frac{i
\theta}{2}}, \qquad
\rho>0, \quad \theta \in [0,2\pi),\\
\label{branch cut} (\rho e^{i \theta})^{1/3} &=& \rho^{1/3}\,
e^{\frac{i \theta}{3}}, \qquad \rho>0,\quad \theta \in (-\pi, +
\pi].
\end{eqnarray}
\end{theorem}
\begin{remark}
In the case that the recurrence coefficients do not depend on $N$
the existence of the limit~\reff{ratiolimit} was already proven in
\cite{Aptekarev1}.  Our proof will be based on similar arguments.
\end{remark}
\begin{remark}
\label{gevolg ratio-as} %
Under the conditions of Theorem~\ref{ratioas}, by taking the
derivative of \reff{ratiolimit} we get
\begin{equation}
\label{lemmaeenbijasymp} \lim_{n/N \to t}  \left(
\frac{P_{n,N}'(z)}{P_{n,N}(z)}
 -\frac{P_{n+1,N}'(z)}{P_{n+1,N}(z)}\right)
 = \left\{
\begin{array}{ll}
{\ds \frac{1}{\alpha} \frac{\phi '(z/\alpha)}{\phi(z/\alpha)}}, &
\qquad \alpha>0,\\[3ex]
{\ds -\frac{1}{z}}, & \qquad \alpha=0,
\end{array}
\right.
\end{equation}
uniformly on compact subsets of $\mathbb{C}\setminus [m,M]$.
\end{remark}
In order to prove Theorem~\ref{ratioas} we need part (a) of the
following lemma. The whole lemma will be used in the proof of
Theorem \ref{as zero distr} as well. It can be found in, e.g.,
\cite[Lemma~2.2]{Kuijlaars1} but we include a short proof for
completeness.

\begin{lemma}
\label{leminterlalg} %
Suppose that the zeros of the monic
polynomials $p_{n-1}$ and $p_n$, with degree $n-1$ and $n$,
respectively, are simple and real, interlace and lie in $[m,M]$. Then %
\begin{itemize}
\item[(a)] $\left|{\ds
                \frac{p_{n-1}(z)}{p_{n}(z)}}\right|  \le {\ds \frac{1}{\dist
                (z,[m,M])}}$,  \qquad  $\forall z\in \mathbb{C}\setminus
                [m,M]$, \\
\item[(b)] $\left|{\ds \frac{p_{n-1}(z)}{p_{n}(z)}}\right| \ge
                {\ds \frac{1}{2|z|}}$,  \qquad
                if  $|z|>\max (|m|,|M|)$.
\end{itemize}
\end{lemma}
\begin{proof}
Denote the real zeros of $p_n$ by $y_1,\ldots,y_n$.  Since
$p_{n-1}$ and $p_n$ are monic and their zeros interlace, there
exist $w_j>0$, $\sum_{j=1}^n w_j=1$, so that
\[\frac{p_{n-1}(z)}{p_{n}(z)}=\sum_{j=1}^n\frac{w_j}{z-y_j}.\]
Then note that, because $y_j\in [m,M]$, for all $z\in
\mathbb{C}\setminus [m,M]$ we have $|z-y_j|\ge \dist (z,[m,M])$,
$1\le j\le n$.  This immediately proves part (a) of the lemma.

If $|z|>\max (|m|,|M|)$, then $|y_j/z|<1$ and therefore $\Re (
\frac{1}{1-y_j/z})>\frac{1}{2}$, $1\le j\le n$.  Hence
\[\frac{1}{|z|}\left|\sum_{j=1}^n\frac{w_j}{1-y_j/z}\right|
\ge \frac{1}{|z|}\,\Re
\left(\sum_{j=1}^n\frac{w_j}{1-y_j/z}\right)
> \frac{1}{2|z|}\sum_{j=1}^n w_j=\frac{1}{2|z|},\]
which proves part (b).
\end{proof}
We also need the following properties of the function $\phi$.
\begin{lemma}
\label{properties phi}
The function $\phi$ is analytic on
$\mathbb{C} \setminus [0,1]$ and satisfies
\begin{itemize}
\item[(a)] $z\phi(z)=\left(1+{\ds \frac{4 \phi(z)}{27}}\right)^3$,
\qquad $z\in \mathbb{C} \setminus [0,1]$, \\%
\item[(b)] $\phi(z)=
z^{-1} +\mathcal{O}(z^{-2})$, \qquad as $z\to \infty$.
\end{itemize}
\end{lemma}

\begin{proof} %
By the choice of branch cuts for the square and the cubic root,
see \reff{branch cut}, the function $\phi$ is certainly analytic
on $\mathbb{C}\setminus (-\infty,1]$.  For $x<0$, a simple
calculation also shows that $\lim_{\epsilon\downarrow
0}\phi(x+i\epsilon)=\lim_{\epsilon\downarrow 0}\phi(x-i\epsilon)$.
So $\phi$ is analytic on $\mathbb{C}\setminus [0,1]$.
If we define
\begin{equation}
\label{uminplus} %
u_{\pm}(z):=-1 \pm \sqrt{1-z}
\end{equation}
then
\begin{equation}
\label{somproduct} u_+(z)u_-(z)=z,\qquad u_+(z)+u_-(z)=-2.
\end{equation}
Using this gives
\begin{eqnarray*}
\left(1+{\ds \frac{4 \phi(z)}{27}}\right)^3 %
&=& \frac{27}{8}\, z \left(\omega_3 (u_+(z))^{1/3} + (u_-(z))^{1/3}\right)^3\\
& & \\
&=& \frac{27}{8}\,z \left(-2 + 3 \omega_3 z^{1/3}\left(\omega_3
(u_+(z))^{1/3}+(u_-(z))^{1/3} \right)\right),
\end{eqnarray*}
which verifies part (a) of the lemma.

Next, we will prove that $\phi$ tends to zero as $z\to \infty$.
From part (a) we then obtain that $\lim_{z\to \infty} z\phi(z)=1$
and, since $\phi$ is analytic in a neighborhood of infinity, this
implies part (b) of the lemma. Applying the formula
$(a+b)(a^2+b^2-ab)=a^3+b^3$ and \reff{somproduct}, we observe that
\begin{equation}
\label{truuk} \frac{4}{27}\,\phi(z)= -\frac{3\,\omega_3\,
z^{1/3}}{\omega_3^2 \left[(u_+(z))^{1/3}\right]^2 +
\left[(u_-(z))^{1/3}\right]^2 - \omega_3 z^{1/3}}-1.
\end{equation}
Now take for a moment $z=1+L$, with $L>0$.  By the definition of
the square root we have $u_{\pm}(1+L)=-1 \pm i\sqrt{L}$. So we can
write
\begin{eqnarray*}
u_+(1+L)&=&\rho(L)\, \e^{i
\left(\frac{\pi}{2}+\varepsilon(L)\right)},\\
u_-(1+L)&=&\rho(L)\, \e^{i
\left(-\frac{\pi}{2}-\varepsilon(L)\right)}, \qquad
0<\varepsilon(L)<\frac{\pi}{2}.
\end{eqnarray*}
Obviously we then get
\begin{eqnarray}
\label{stuk1} %
\omega_3^2 \left[(u_+(1+L))^{1/3}\right]^2 &=& (\rho(L))^{2/3}\,
\e^{i\left(-\frac{\pi}{3}+\frac{2}{3}\, \varepsilon(L)\right)},\\
\label{stuk2} %
\left[(u_-(1+L))^{1/3}\right]^2 &=&  (\rho(L))^{2/3}\,
\e^{i\left(-\frac{\pi}{3}-\frac{2}{3}\, \varepsilon(L)\right)}.
\end{eqnarray}
Finally, notice that
\begin{equation}
\lim_{L\to +\infty} \varepsilon(L)=0 \quad \mbox{and}\quad \rho(L)
\sim \sqrt{L}, \quad L \to +\infty.
\end{equation}
Hence, combining \reff{truuk}, \reff{stuk1} and \reff{stuk2}, we
obtain
\[\lim _{L \to +\infty} \frac{4}{27}\,\phi(1+L)=\lim_{L \to +\infty}
\frac{3}{\e^{i\left(\frac{2}{3}\, \varepsilon(L)\right)}
+\e^{-i\left(\frac{2}{3}\, \varepsilon(L)\right)}+1}-1 =0.\]
Since
$\phi$ is analytic in a neighborhood of infinity, then also
$\lim_{z\to \infty} \phi(z)=0$.
\end{proof}
Now we give the proof of Theorem \ref{ratioas}.
\begin{varproof}
\textbf{of Theorem~\ref{ratioas}:}
It is enough to prove the cases $\alpha=0$ and $\alpha=1$. The
more general case $\alpha>0$ is then obtained by taking $\tilde
P_{n,N}(z):=P_{n,N}\left(\alpha z\right)/\alpha^n$.

We first prove the case $\alpha=1$. By the assumptions on the
zeros of the polynomials $P_{n,N}$ every member of
\begin{equation}
\label{normalefamilie} %
\left\{ \left.\frac{P_{n,N}(z)}{P_{n+1,N}(z)} \ \right|\  n,N \in
\mathbb{N},\,  n \leq t^{*} N \right\}
\end{equation}
satisfies the estimate in part (a) of Lemma~\ref{leminterlalg}.
So, the family \reff{normalefamilie} is uniformly bounded on
compact subsets of $\overline{\mathbb{C}} \setminus [m,M]$.  By
the theorem of Montel \cite[p.~563]{henrici} we then know that
\reff{normalefamilie} is a normal family on $\overline{\mathbb{C}}
\setminus [m,M]$.  For a sequence $\{(n_j,N_j)\}_{j\ge 1}$, with
$n_j,N_j \to \infty$, $n_j/N_j\to t$  as $j\to \infty$, we have
that, if $j$ is sufficiently large, the function
\begin{equation}
\label{f} %
f_j(z):=\frac{P_{n_j,N_j}(z)}{P_{n_j+1,N_j}(z)}
\end{equation}
belongs to the normal family (\ref{normalefamilie}).  The
corresponding sequence $\{f_j\}_{j\ge 1}$ then has a subsequence
that converges uniformly on compact subsets of
$\overline{\mathbb{C}}\setminus [m,M]$. If we can prove that the
limit of any such subsequence is $\phi$, then, by a standard
compactness argument, the full sequence $\{f_j\}_{j\ge 1}$
converges uniformly on compact subsets of
$\overline{\mathbb{C}}\setminus [m,M]$ to $\phi$.  This then
proves the theorem in the case $\alpha=1$.

We will show that for each sequence $n_i, N_i \to \infty$ with
$n_i/N_i \to t$ such that the functions $\{f_i\}_{i\ge 1}$
converge uniformly on compact subsets of $\overline{\mathbb{C}}
\setminus [m,M]$, we have
\begin{equation}
\label{klimietgedrag} %
f(z):= \lim_{i \to \infty} f_i(z) = \phi(z) + \mathcal{O}(z^{-k}),
\qquad \mbox{as } z \to \infty,
\end{equation}
for each $k\in\mathbb{N}$. The uniqueness of the Laurent expansion
around infinity then implies that $f(z)=\phi(z)$. We show this by
induction on $k$. The case $k=1$ follows from
Lemma~\ref{properties phi}~(b) and $f_i(z) = \mathcal{O}(z^{-1})$,
for every $i\ge 1$. Next, suppose that the claim holds for some $k
\geq 1$ and consider a sequence $\{(n_i,N_i)\}_{i\ge 1}$ such that
$n_i,N_i \to \infty$, $n_i/N_i\to x$ and the functions
$\{f_i\}_{i\ge 1}$ converge uniformly on compact subsets of
$\overline{\mathbb{C}}\setminus [m,M]$ to some function $f$ as
$i\to \infty$.  If we put
\begin{eqnarray*}
g_i(z) &:=& \frac{P_{n_{i}-1, N_i}(z)}{P_{n_{i},N_i}(z)}, \qquad z
\in \mathbb{C} \setminus [m,M],\\
h_i(z) &:=& \frac{P_{n_{i}-2, N_i}(z)}{P_{n_{i}-1,N_i}(z)}, \qquad
z \in \mathbb{C} \setminus [m,M],
\end{eqnarray*}
then from the recurrence relation \reff{recursie ratio as} we
obtain
\begin{equation}
\label{recursie2} %
z= f_i(z)^{-1}+b_{n_i, N_i} + c_{n_i, N_i}\, g_i(z) + d_{n_i,
N_i}\, h_i(z).
\end{equation}
Since $t<t^{*}$ we may assume without loss of generality that $n_i
< t^{*} N_i$ for every $i\ge 1$. Then $\{g_i\}_{i\ge 1}$ and
$\{h_i\}_{i\ge 1}$ are subsets of the normal family
(\ref{normalefamilie}).  Therefore, there is a sequence $i_j \to
\infty$, $j\to \infty$, such that $\{g_{i_j}\}_{j\ge 1}$ and
$\{h_{i_j}\}_{j\ge 1}$ converge uniformly on compact subsets of
$\overline{\mathbb{C}} \setminus [m,M]$ with limit $g$ and $h$,
respectively. If we pass to such a subsequence and take limits in
\reff{recursie2}, then by (\ref{limbcd}) we find
\begin{equation}
\label{fgh} %
z=\frac{1}{f(z)} + 3\left(\frac{4}{27}\right) +
3\left(\frac{4}{27}\right)^2 g(z) + \left(\frac{4}{27}\right)^3
g(z) h(z), \quad z \in \mathbb{C} \setminus [m,M].
\end{equation}
By the induction hypothesis we now have that
\begin{eqnarray*}
g(z) & = & \phi(z) + \mathcal{O}(z^{-k}), \qquad z \to \infty,\\
h(z) & = & \phi(z) + \mathcal{O}(z^{-k}), \qquad z \to \infty.
\end{eqnarray*}
Applying this to \reff{fgh}, by Lemma~\ref{properties phi} we then
get
\begin{equation*}
\frac{1}{f(z)}= z- 3\left(\frac{4}{27}\right) -
3\left(\frac{4}{27}\right)^2 \phi(z) - \left(\frac{4}{27}\right)^3
\phi(z)^2 + \mathcal{O}(z^{-k})=\frac{1}{\phi(z)} +
\mathcal{O}(z^{-k}).
\end{equation*}
Since $\phi(z)=\mathcal{O}(z^{-1})$, this implies
\begin{equation*}
f(z) = \frac{\phi(z)}{1+ \phi(z)\mathcal{O}(z^{-k})}
=\frac{\phi(z)}{1+ \mathcal{O}(z^{-k-1})} = \phi(z) +
\mathcal{O}(z^{-k-2}).
\end{equation*}
So we proved that (\ref{klimietgedrag}) also holds with $k$
replaced by $k+2$. Therefore, it holds for all $k$.

Finally, for the case $\alpha=0$ the proof is similar.  In fact,
\reff{ratiolimit} then easily follows by taking limits in
\reff{recursie2}.
\end{varproof}

\section{Proof of Theorem~\ref{as zero distr}}
\label{proof main theorem}

In order to prove the asymptotic result in Theorem~\ref{as zero
distr} we first have a closer look at the function
$\frac{\phi'}{\phi}$, which is analytic on $\mathbb{C}\setminus
[0,1]$. Here we will use the relation
\begin{equation}
\label{verhouding naar phi} %
\frac{\phi '(z)}{\phi(z)} = \frac{\frac{4\phi(z)}{27}+1}{
z\left(\frac{8\phi(z)}{27}-1\right)} =\frac{1}{2z} +
 \frac{3}{2z}\left(\frac{8\phi(z)}{27}-1\right)^{-1},
\end{equation}
which can be obtained by differentiating part (a) of
Lemma~\ref{properties phi}.
First of all, we are interested in the jump across its branch cut.
\begin{lemma}
\label{lemma jump} %
The jump of the function $\frac{\phi'}{\phi}$
across its branch cut is given by
\begin{equation}
\label{functie m} %
m(x):=\lim_{\varepsilon \downarrow 0} \frac{\phi
'(x+i\epsilon)}{\phi(x+i\varepsilon)} - \lim_{\varepsilon
\downarrow 0} \frac{\phi '(x-i\varepsilon)}{\phi(x-i\epsilon)} =2
\pi i \ \frac{\d \upsilon_{[0,1]}}{\d x}(x), \qquad x \in (0,1),
\end{equation}
where $\upsilon_{[0,1]}$ is defined as in Theorem~\ref{as zero
distr}.
\end{lemma}

\begin{proof}
Let $x \in (0,1)$.  By \reff{verhouding naar phi} we easily obtain
\begin{equation}
\label{m naar phi} %
m(x)= \frac{3}{2 x}  \left\{ \lim_{\varepsilon \downarrow
0}\left(\frac{8\phi(x\e^{i \varepsilon})}{27}-1\right)^{-1}-
\lim_{\varepsilon \downarrow 0}\left(\frac{8\phi(x\e^{-i
\varepsilon})}{27}-1\right)^{-1} \right\}.
\end{equation}
Applying the definitions \reff{branch cut} we get
\begin{eqnarray*}
\lim_{\varepsilon \downarrow 0} \left( -1+ \sqrt{1-x\e^{\pm i
\varepsilon}} \right)^{1/3}
&=&\e^{i \frac{\pi}{3}} \left( 1\pm \sqrt{1-x}\right)^{1/3},\\
\lim_{\varepsilon \downarrow 0} \left( -1- \sqrt{1-x\e^{\pm i
\varepsilon}} \right)^{1/3} &=&\e^{-i \frac{\pi}{3}} \left( 1 \mp
\sqrt{1-x}\right)^{1/3}.
\end{eqnarray*}
Using the notation $v_{\pm}(x):=(1 \pm \sqrt{1-x})^{1/3}$ we then
have
\[\lim_{\varepsilon \downarrow 0}\frac{8\phi(x\e^{\pm i
\varepsilon})}{27}-1=3x^{1/3}\left\{\e^{-\frac{i\pi}{3}}v_{\pm}(x)
+e^{\frac{i\pi}{3}}v_{\mp}(x)\right\}-3.\]
So, also applying the relations $x^{1/3}=v_+(x)v_-(x)$ and
$2=v_+(x)^3+v_-(x)^3$, equation~\reff{m naar phi} becomes
\begin{eqnarray*}
m(x) & = & \frac{2\sqrt{3}i}{x^{2/3}}
\frac{v_+(x)-v_-(x)}{\left(x^{1/3}(v_+(x)+v_-(x))-2\right)^2+3x^{2/3}\left(v_+(x)-v_-(x)\right)^2}\\
 & = & \frac{2\sqrt{3}i}{x^{2/3}}\, (v_+(x)-v_-(x))^{-1}
 \left[(v_-(x)^2-v_+(x)^2)^2 +3v_+(x)^2v_-(x)^2\right]^{-1}.
\end{eqnarray*}
If we multiply the numerator and denominator both by
$v_+(x)+v_-(x)$, then we finally obtain
\[m(x) =  \frac{2\sqrt{3}i}{x^{2/3}}
\frac{v_+(x)+v_-(x)}{\left(v_+(x)^3+v_-(x)^3\right)\left(v_+(x)^3-v_-(x)^3\right)}
= \frac{\sqrt{3}i}{2x^{2/3}} \frac{v_+(x)+v_-(x)}{\sqrt{1-x}},\]
which proves \reff{functie m}.
\end{proof}
A second point of interest is the behavior of the function
$\frac{\phi'}{\phi}$ in the neighborhood of its branch points $0$
and $1$.

\begin{lemma}
\label{branch points} %
Near the points $0$ and $1$ we have
\begin{eqnarray}
\label{bijnul} \left|\frac{\phi'(\varepsilon
\e^{i\theta})}{\phi(\varepsilon \e^{i\theta})}\right| & = &
\mathcal{O}\left(\varepsilon^{-2/3}\right), \quad \varepsilon
\downarrow 0, \qquad \theta\in
\left(0,2\pi \right),\\
\label{bij1} \left|\frac{\phi'(1+\varepsilon
\e^{i\theta})}{\phi(1+\varepsilon \e^{i\theta})}\right| & = &
\mathcal{O}\left(\varepsilon^{-1/2}\right), \quad \varepsilon
\downarrow 0, \qquad \theta \in \left(-\pi,\pi \right).
\end{eqnarray}
\end{lemma}

\begin{proof}
For $\theta \in \left(0,2\pi \right)$ we easily see that
\[
\lim_{\varepsilon \downarrow 0} \phi(\varepsilon \e^{i\theta}) =
-\frac{27}{4},\qquad \left|\frac{4 \phi(\varepsilon
\e^{i\theta})}{27}+1\right| = \mathcal{O}\left(\varepsilon^{1/3}
\right), \quad \varepsilon \downarrow 0.
\]
Applying this to the first equality in \reff{verhouding naar phi}
then gives expression (\ref{bijnul}).

We now take $\theta \in \left(-\pi,\pi \right)$. We then have $
u_\pm (1+\varepsilon e^{i\theta}) = -1 \pm \sqrt{ \varepsilon}
\,\e^{i \frac{\pi+\theta}{2}}$, where we use the notation
\reff{uminplus}, and applying the definition of the third root
\reff{branch cut} we get
\begin{equation}
\label{deel 1} %
\lim_{\varepsilon \downarrow 0} \phi(1+\varepsilon
\e^{i\theta}) =
\frac{27}{4}\left(\frac{3\omega_3}{2}\left\{\omega_3\e^{\frac{\pi
i}{3}}+\e^{-\frac{\pi i}{3}}\right\}-1\right)=\frac{27}{8}.
\end{equation}
If we write $ u_\pm (1+\varepsilon e^{i\theta}):= \rho_\pm
(\varepsilon,\theta) \e^{i\eta_\pm (\varepsilon,\theta)} $,
meaning the polar coordinates, then a closer look gives
\[\eta_\pm (\varepsilon,\theta) = \pm \pi \mp \sqrt{\varepsilon} \cos \left(\theta/2\right)
+\mathcal{O}(\varepsilon) , \quad \varepsilon \downarrow 0,\]
and
\begin{align}
& \e^{i \frac{\eta_\pm (\varepsilon,\theta)}{3}}  =
\e^{\pm\frac{\pi i}{3}} \left(1 \mp \frac{i\sqrt{\varepsilon}}{3}
\cos \left(\theta/2\right)\right)+\mathcal{O}(\varepsilon)
, \quad \varepsilon \downarrow 0,\\[1ex]
& \rho_\pm (\varepsilon, \theta)^{1/3}  =  1 \pm
\frac{\sqrt{\varepsilon}}{3} \sin
\left(\theta/2\right)+\mathcal{O}(\varepsilon) , \quad \varepsilon
\downarrow 0.
\end{align}
From this we easily obtain
\begin{equation}
\label{deel 2} %
\frac{8 \phi(1+\varepsilon \e^{i\theta})}{27}-1 =
-\sqrt{3\varepsilon}
\e^{\frac{i\theta}{2}}+\mathcal{O}(\varepsilon) , \quad
\varepsilon \downarrow 0.
\end{equation}
Applying \reff{deel 1} and \reff{deel 2} to the first equality in
\reff{verhouding naar phi} then finally leads to (\ref{bij1}).
\end{proof}
As a corollary of Lemma~\ref{lemma jump} and Lemma~\ref{branch
points} we obtain that $\frac{\phi'}{\phi}$ is the Stieltjes
transform of the measure $\upsilon_{[0,1]}$, up to a minus sign.
\begin{lemma}
\label{stieltjes} %
Let $\phi$ be defined by
(\ref{phidefinitie}), then
\begin{equation}
\frac{\phi '(z)}{\phi(z)} = -\int \frac{1}{z-x} \, \d
\upsilon_{[0,1]}(x), \qquad z \in \mathbb{C}\setminus [0,1],
\end{equation}
with $\upsilon_{[0,1]}$ defined as in Theorem~\ref{as zero distr}.
\end{lemma}

\begin{proof}
By Lemma~\ref{lemma jump} and Lemma~\ref{branch points} and
applying Lemma~\ref{properties phi}~(b) to \reff{verhouding naar
phi}, the function $\frac{\phi'}{\phi}$ satisfies the following
additive Riemann-Hilbert problem:\\[-1ex]
\begin{itemize}
\item[(P1)] $f$ is analytic in $\mathbb{C}\setminus [0,1]$, \\%
\item[(P2)] $\lim\limits_{\varepsilon \downarrow 0}
f(x+i\varepsilon)-\lim\limits_{\varepsilon \downarrow 0}
f(x-i\varepsilon)=2\pi i \ \frac{\d \upsilon_{[0,1]}}{\d x}(x)$,
\quad for $x \in (0,1)$, \\%
\item[(P3)]
$f(z)=-z^{-1}+\mathcal{O}\left(z^{-2}\right)$, \quad as $z\to
\infty$, \\%
\item[(P4)] $f(z)=\mathcal{O}\left(z^{-2/3}\right)$, \quad as
$z\to 0$, \\ $f(z)=\mathcal{O}\left((1-z)^{-1/2}\right)$, \quad as
$z\to 1$.
\end{itemize}
\ \\[-1ex]
If $f$ and $g$ are both solutions of this Riemann-Hilbert
problem, then it is easily seen that $f-g$ is analytic in
$\mathbb{C}\setminus \{0,1\}$. Moreover, $0$ and $1$ are removable
singularities by (P4). Liouville's Theorem and (P3) then imply
$f\equiv g$, meaning that the Riemann-Hilbert problem has a unique
solution. So it is enough to show that
    \[f(z):=\int \limits \frac{1}{x-z} \,\d \upsilon_{[0,1]}(x),\qquad z\in\mathbb{C}\setminus [0,1],\]
satisfies (P1)-(P4). Properties (P1) and (P3) easily follow from
the fact that $\frac{\d \upsilon_{[0,1]}}{\d x}$ is a probability
measure on $(0,1)$. By the Plemelj-Sokhotskii formula for Cauchy
integrals, see, e.g., \cite[p.~43, (18.1)]{Musk}, $f$ satisfies
(P2). Finally, the behavior at the branch points, see (P4), easily
follows from \cite[p.~74, (29.5) and (29.6)]{Musk} and
\reff{upsilon}.
\end{proof}
\begin{remark}
\label{stieltjes algemeen} %
As an easy consequence of
Lemma~\ref{stieltjes} we obtain
\begin{equation}
-\int \frac{1}{z-x} \, \d\upsilon_{[0,\alpha]}(x) %
= \left\{
\begin{array}{ll}
{\ds \frac{1}{\alpha}\frac{\phi '(z/\alpha)}{\phi(z/\alpha)}}, &
\quad \alpha>0,\\[2ex]
-{\ds \frac{1}{z}}, & \quad \alpha=0,
\end{array}
\right.
\qquad z \in \mathbb{C}\setminus [0,\alpha],
\end{equation}
where $\upsilon_{[0,\alpha]}$ is defined as in Theorem~\ref{as
zero distr}.
\end{remark}
We are now ready to prove Theorem~\ref{as zero distr}.
\begin{varproof}
\textbf{of Theorem~\ref{as zero distr}:} Let $t>0$ and fix a
number $t^\star>t$. Clearly, the convergence (\ref{limbcd met t})
and the fact that the function $\alpha$ is continuous on
$[0,\infty)$ imply that the recurrence coefficients are uniformly
bounded if $n/N$ is restricted to compact subsets of $[0,\infty)$.
So,
\begin{equation}
0<R:=\sup\{1+|b_{n,N}|+|c_{n,N}|+|d_{n,N}|\, :\, n\le t^\star
N\}<+\infty.
\end{equation}
By the recurrence \reff{four term recursie} we have
$P_{n,N}(z)=\det(zI_n-L_{n,N})$, with
\begin{equation}
\label{matrix L} %
L_{n,N}  =  \left(
\begin{array}{cccccc}
b_{0,N}&1 & 0 & \ldots & \ldots &  0\\
c_{1,N} &b_{1,N} &1 &\ddots &  &\vdots \\
d_{2,N} & c_{2,N} & b_{2,N} &1 & \ddots  & \vdots\\
0& d_{3,N} & c_{3,N} & b_{3,N}  & \ddots & 0\\
\vdots&\ddots &\ddots &\ddots  &\ddots &1 \\
0&\ldots &0 & d_{n-1,N} & c_{n-1,N} & b_{n-1,N}
\end{array} \right).
\end{equation}
As a consequence,  the zeros of $P_{n,N}$ are bounded by
$\|L_n\|_\infty$.  For $n\le t^\star N$ we then have that the
simple zeros of $P_{n,N}$ lie in the interval $[-R,R]$. Moreover,
they are assumed to satisfy the interlacing property (for fixed
$N$) and one can observe that
\[\alpha(st)\le
\left(1+\frac{4\alpha(st)}{27}\right)^3\le R, \qquad \mbox{for
each } 0<s\le 1.\]
So, by Remark~\ref{gevolg ratio-as} and
Remark~\ref{stieltjes algemeen} we establish
\begin{equation}
\label{Ratio with s} %
\lim_{n/N \to t}  \left( \frac{P_{\lfloor sn
\rfloor+1,N}'(z)}{P_{\lfloor sn \rfloor+1,N}(z)} -\frac{P_{\lfloor
sn \rfloor,N}'(z)}{P_{\lfloor sn \rfloor,N}(z)}\right) = \int
\frac{1}{z-x} \, \d\upsilon_{[0,\alpha(st)]}(x)
\end{equation}
uniformly on compact subsets of $ \mathbb{C}\setminus [-R,R]$,
where $0<s\le 1$ and $\lfloor sn \rfloor$ denotes the greatest
integer less than or equal to $sn$.

Note that
\begin{equation}
\label{som naar integraal} %
\frac{1}{n} \frac{P_n'(z)}{P_n(z)} = \frac{1}{n} \sum_{k=0}^{n-1}
\left( \frac{P_{k+1}'(z)}{P_{k+1}(z)} -\frac{P_k'(z)}{P_k(z)}
\right) =  \int_{0}^{1 }\left( \frac{P_{\lfloor sn
\rfloor+1}'(z)}{P_{\lfloor sn \rfloor+1}(z)} -\frac{P_{\lfloor sn
\rfloor}'(z)}{P_{\lfloor sn \rfloor}(z)} \right) \d s.
\end{equation}
For $n\le t^\star N$ the zeros of the polynomials $P_{n,N}$ are
simple, lie in $[-R,R]$ and satisfy the interlacing property for
fixed $N$.  From Lemma~\ref{leminterlalg}~(b) we then get
\[\left|\frac{P_{\lfloor sn
\rfloor+1,N}(z)}{P_{\lfloor sn \rfloor,N}(z)}\right| \le
2|z|,\qquad |z|>R.\]
With a similar argument as in Lemma~\ref{leminterlalg}~(a) we can
also prove
\[\left|\left(\frac{P_{\lfloor sn
\rfloor,N}(z)}{P_{\lfloor sn \rfloor+1,N}(z)}\right)'\right| \le
\frac{1}{\dist (z,[-R,R])^2}, \qquad z\in\mathbb{C}\setminus
[-R,R].\]
Combining these two results we have, for $|z|>R$,
\begin{eqnarray}
\nonumber \left|\frac{P_{\lfloor sn \rfloor+1}'(z)}{P_{\lfloor sn
\rfloor+1}(z)} -\frac{P_{\lfloor sn \rfloor}'(z)}{P_{\lfloor sn
\rfloor}(z)}\right| %
& = & %
\left|\frac{P_{\lfloor sn \rfloor+1,N}(z)}{P_{\lfloor sn
\rfloor,N}(z)}\right| \left|\left(\frac{P_{\lfloor sn
\rfloor,N}(z)}{P_{\lfloor sn
\rfloor+1,N}(z)}\right)'\right| \\ %
& \le & %
\frac{2|z|}{\dist (z,[-R,R])^2}.
\end{eqnarray}
So, we can apply Lebesgue's dominated convergence theorem on
(\ref{som naar integraal}) and by (\ref{Ratio with s}) we obtain
\begin{eqnarray}
\nonumber \lim_{n/N\to t} \int \frac{1}{z-x} \, \d\nu(P_{n,N})(x)
& = & \int_0^1 \int \frac{1}{z-x} \, \d\upsilon_{[0,\alpha(st)]}(x) \,\d s\\
\label{finale gelijkheid stieltjes} & = & \frac{1}{t} \int_0^t
\int \frac{1}{z-x} \, \d\upsilon_{[0,\alpha(s)]}(x) \d s,
\end{eqnarray}
for $|z|>R$. By \cite[Theorem~2]{Geronimo}, which is a gloss on
the theorem of Grommer and Hamburger \cite[p.~104-105]{Wintner},
we then finally establish \reff{limitnu}.
\end{varproof}

\end{document}